\documentclass[10pt]{article}
\usepackage{amstext}
\usepackage{amsfonts}
\usepackage{amssymb}
\usepackage{amsbsy}
\usepackage{amsmath}
\usepackage{latexsym}
\usepackage{xy}
\usepackage{hhline}
\xyoption{all}
\newcommand{\vtx}[1]{*+[o][F-]{\scriptscriptstyle #1}} 

\newcounter{num}[section] %

\newenvironment{theo}
{\refstepcounter{num}%
\bigskip\noindent{\bf Theorem~\arabic{section}.\arabic{num}. }\it}

\newenvironment{cor}
{\refstepcounter{num}%
\bigskip\noindent{\bf Corollary~\arabic{section}.\arabic{num}. }\it}

\newenvironment{lemma}
{\refstepcounter{num}%
\bigskip\noindent{\bf Lemma~\arabic{section}.\arabic{num}. }\it}

\newenvironment{example}
{\refstepcounter{num}%
\bigskip\noindent{\bf Example~\arabic{section}.\arabic{num}.}}

\newenvironment{remark}
{\refstepcounter{num}%
\bigskip\noindent{\bf Remark~\arabic{section}.\arabic{num}.}\it}

\newcommand{\Ref}[1]{(\ref{#1})}

\newcounter{thepic}

\newenvironment{proof}{\medskip\noindent{\it Proof. }}
{$\Box$ \bigskip}
\newenvironment{proof_lemma_reduction}
{\noindent\textbf{Proof of Lemma~\ref{lemma_reduction}. }} {$\Box$ \bigskip}
\newenvironment{max_path_problem}{\bigskip\noindent{\bf Guide's Problem. }}
{\bigskip}
\newenvironment{eq}{\begin{equation}}{\end{equation}}

\newcommand{\si}{\sigma}
\newcommand{\al}{\alpha}

\newcommand{\la}{\lambda}
\newcommand{\de}{\delta}

\newcommand{\ch}[1]{\check{#1}}
\newcommand{\ov}[1]{\overline{#1}}
\newcommand{\un}[1]{{\underline{#1}} }

\newcommand{\tr}{\mathop{\rm tr}}

\newcommand{\mdeg}{\mathop{\rm mdeg}}
\newcommand{\degII}[2]{\mathop{\rm {deg}}^{o}_{#1}(#2)}
\newcommand{\diag}{\mathop{\rm diag}}
\newcommand{\Char}{\mathop{\rm char}}
\newcommand{\sign}{\mathop{\rm{sgn }}}

\newcommand{\NN}{{\mathbb{N}} }
\newcommand{\ZZ}{{\mathbb{Z}} }

\newcommand{\Q}{\mathcal{Q}}    
\newcommand{\QG}{\mathcal{G}}    
\newcommand{\n}{\boldsymbol{n}} 
\newcommand{\Ver}[1]{\mathop{{\rm ver}(#1)}} 
\newcommand{\Arr}[1]{\mathop{{\rm arr}(#1)}} 
\newcommand{\Inv}{I}       
\newcommand{\path}{\mathop{\rm path}} 

\newcommand{\QuiverNull}{\Q}     

\newcommand{\Symm}{\mathcal{S}}                 
\newcommand{\loopR}[3]{%
\begin{picture}(20,0)(#1,#2)
\put(-2,1){\llap{$\scriptstyle #3$}} \put(10,3){\circle{20}} \put(20,6){\vector(1,-4){1}}
\end{picture}}
\newcommand{\loopL}[3]{%
\begin{picture}(20,0)(#1,#2)
\put(22,1){$\scriptstyle #3$} \put(10,3){\circle{20}} \put(0,6){\vector(-1,-4){1}}
\end{picture}}
\newcommand{\loopD}[3]{%
\begin{picture}(20,20)(#1,#2)
\put(22,10){$\scriptstyle #3$} \put(10,10){\circle{20}} \put(10,0){\vector(-4,-1){1}}
\end{picture}}

\newcommand{\rombL}{
\begin{picture}(0,0)
\put(80,0){\vector(-4,3){40} \put(13,15){$\scriptstyle 2$}}%
\put(40,30){\vector(0,-1){22}}\put(40,-8){\vector(0,-1){22}}%
\put(40,-30){\vector(-4,3){40}\put(3,10){$\scriptstyle t$}}%
\put(40,-30){\vector(4,3){40} \put(-17,10){$\scriptstyle 2$}}%
\put(40,4){\circle*{1}}\put(40,0){\circle*{1}}\put(40,-4){\circle*{1}}%
\put(-10,31){
{\xymatrix@C=1.3cm@R=0.9cm{ %
& \ar@/^/@{<-}[ld] \ar@/_/@{<-}[ld] \\%
\\}}}%
\put(19,16){\put(2,-2){\circle*{1}}\put(0,0){\circle*{1}}\put(-2,2){\circle*{1}}}%
\put(11,22){$\scriptstyle 1$}\put(25,4){$\scriptstyle 1$}\put(-6,-2){$\scriptstyle u$}\put(39,32){$\scriptstyle v$}%
\end{picture}} %

\newcommand{\cycleSmall}{
\begin{picture}(0,0)
\put(-5,0){\xymatrix@C=1.75cm@R=2cm{ %
\ar@/^/@{->}[r]&\\
\\
}}%
\put(50,0){\vector(-1,-2){10}}%
\put(10,-20){\vector(-1,2){10}}%
\put(25,-20){\put(0,0){\circle*{1}}\put(4,0){\circle*{1}}\put(-4,0){\circle*{1}}}%
\end{picture}} %

\begin{document}
\title{Indecomposable invariants of quivers for dimension $(2,\ldots,2)$ and maximal paths.}
 \author{
A.A. Lopatin \\
Institute of Mathematics, SBRAS, \\
Pevtsova street, 13,\\
Omsk 644099 Russia \\
artem\underline{ }lopatin@yahoo.com \\
http://www.iitam.omsk.net.ru/\~{}lopatin/\\
}
\date{} 
\maketitle

\begin{abstract}
An upper bound on degrees of elements of a minimal generating system for invariants of
quivers of dimension $(2,\ldots,2)$ is established over a field of arbitrary
characteristic and its precision is estimated. The proof is based on the reduction to the
problem of description of maximal paths satisfying certain condition.
\end{abstract}




\section{Introduction}\label{section_intro}
We work over an infinite field $K$ of arbitrary characteristic $\Char(K)$. All vector
spaces, algebras, and modules are over $K$ unless otherwise stated and all algebras are
associative.

A {\it quiver} $\Q=(\Ver{\Q},\Arr{\Q})$ is a finite oriented graph, where $\Ver{\Q}$ is
the set of vertices and $\Arr{\Q}$ is the set of arrows. For an arrow $a$ denote by $a'$
its head and denote by $a''$ its tail. Loops and multiple arrows are
allowed. The notion of quiver was introduced by Gabriel in~\cite{Gabriel_1972} as an
effective mean for description of different problems of the linear algebra.

For a quiver $\Q$ and a {\it dimension vector} $\n=(\n_{v}\,|\,v\in\Ver{\Q})$ denote by
$\Inv(\Q,\n)$ the algebra of invariants of representations of $\Q$. Invariants of quivers
are important not only in the invariant theory but also in the representational theory
because these invariants distinguish semi-simple representations of a quiver. The algebra
$\Inv(\Q,\n)$ is embedded into the algebra of (commutative) polynomials
$K[x_{ij}(a)\,|\,a\in\Arr{\Q},\ 1\leq i\leq \n_{a'},\,1\leq j\leq \n_{a''}]$. Denote by %
$$X_{a}=\left(
\begin{array}{ccc}
x_{1,1}(a)&\cdots &x_{1,n_{a''}}(a)\\
\vdots&&\vdots\\
x_{n_{a'},1}(a)&\cdots &x_{n_{a'},n_{a''}}(a)\\
\end{array}
\right)$$ %
the $n_{a'}\times n_{a''}$ {\it generic} matrix and by $\si_k(X)$ the $k$-th coefficient
in the characteristic polynomial of an $n\times n$ matrix $X$,
i.e., %
$$\det(\la E-X)=\la^n-\si_1(X)\la^{n-1}+\cdots+(-1)^n\si_n(X).
$$%
In particular, $\si_1(X)=\tr(X)$ and $\si_n(X)=\det(X)$.

For a real number $\alpha$ let $[\alpha]$ be the greatest integer that does not exceed
$\alpha$. We write $\de(i,j)$ for the Kronecker symbol and $\#S$ for the cardinality of a
set $S$.

Let us recall that $a=a_1\cdots a_s$ is a {\it path} in $\Q$ (where
$a_1,\ldots,a_s\in\Arr{\Q}$), if $a_1'=a_2'',\ldots, a_{s-1}'=a_s''$; and $a$ is a {\it
closed} path in a vertex $v$, if $a$ is a path and $a_1''=a_s'=v$. The head of the path
$a$ is $a'=a_s'$ and the tail of $a$ is $a''=a_1''$. Denote
$\Ver{a}=\{a_1'',a_1',\ldots,a_s'\}$, $\Arr{a}=\{a_1,\ldots,a_s\}$, and $\deg(a)=s$. If
$a$ is a closed path, then define the degree of $a$ in a vertex $w$ by
$\deg_w(a)=\#\{i\,|\,a_i'=w,\,1\leq i\leq s\}$. A closed path $a$ is called {\it
primitive} if $\deg_w(a)=1$ for all $w\in\Ver{a}$, i.e., $a$ is without
self-intersections. Denote by $m(\Q)$ the maximal degree of primitive closed paths in
$\Q$. Closed paths $a_1,\ldots,a_s$ in $Q$ are called {\it incident} if
$a_1'=\cdots=a_s'$.

In~\cite{Donkin_1994} Donkin proved that $K$-algebra $\Inv(\Q,\n)$ is generated by
$\si_k(X_{a_s}\cdots X_{a_1})$ for all closed paths $a=a_1\cdots a_s$ in $\Q$ (where
$a_1,\ldots,a_s\in\Arr{\Q}$) and $1\leq k\leq n_{a'}$. For a field of characteristic zero
generators for $\Inv(\Q,\n)$ were described earlier by Le Bruyn and Procesi
in~\cite{Le_Bruyn_Procesi_1990}. Relations between generators are described by Zubkov's  Theorem (see~\cite{Zubkov_Fund_Math_2001}), which for a quiver with one vertex and $\Char(K)=0$ was
independently proven by Razmyslov in~\cite{Razmyslov_1974} and Procesi
in~\cite{Procesi_1976}; for an arbitrary quiver and a field of characteristic zero it was proven by Domokos in~\cite{Domokos_1998}. Notice that $\Inv(\Q,\n)$ has a
grading by degrees that is given by the formula: $\deg(\si_k(X_{a_s}\cdots X_{a_1}))=ks$.

By the Hilbert--Nagata Theorem on invariants, $\Inv(\Q,\n)$ is a finitely generated
graded algebra. But the mentioned generating system is not finite. So it gives rise to
the problem to find out a minimal (by inclusion) homogeneous system of generators
(m.h.s.g.). Given an $\NN$-graded algebra $A$, where $\NN$ stands for non-negative
integers, denote by $A^{+}$ the subalgebra generated by elements of $A$ of positive
degree. It is easy to see that a set $\{a_i\} \subseteq A$ is a m.h.s.g.~if and only if
$\{\ov{a_i}\}$ is a basis of $\ov{A}={A}/{(A^{+})^2}$. An element $a\in A$ is called {\it
decomposable} if it belongs to the ideal $(A^{+})^2$. In other words, a decomposable
element is equal to a polynomial in elements of strictly lower degree. Therefore the
least upper bound $D(\Q,\n)$ for the degrees of elements of a m.h.s.g.~of $\Inv(\Q,\n)$ is equal to the highest degree of indecomposable invariants. In this paper we establish an upper bound on $D(\Q,\n)$ for an arbitrary quiver $\Q$ and $\n=(2,2,\ldots,2)$
and estimate its precision.

In characteristic zero case $\Inv(\Q,\n)$ is generated by invariants of degree at most
$(\sum\n_a)^2$, where the sum ranges over all $a\in\Arr{\Q}$
(see~\cite{Le_Bruyn_Procesi_1990}). All the rest of known results on finite generating
systems for $\Inv(\Q,\n)$ concern a quiver $\Q$ with one vertex and several loops. Some
of these results are presented in Section~\ref{section_matrix_invariants}.

If $\Q_1$ and $\Q_2$ are quivers with $\Ver{\Q_1}\subset \Ver{\Q_2}$ and
$\Arr{\Q_1}\subset \Arr{\Q_2}$, then we say that $\Q_1$ is a {\it subquiver} of $\Q_2$ and
write $\Q_1\subset \Q_2$. A quiver $\Q$ is said to be {\it strongly connected} if there
exists a closed path in $\Q$ that contains all vertices of $\Q$. A quiver with one vertex
and no arrows is also called strongly connected. For a quiver $\Q$ let $\Q_1,\ldots,\Q_k$
be its {\it strongly connected components}, i.e., $\Q_1,\ldots,\Q_k$ are strongly 
connected subquivers of $\Q$, $\Ver{\Q}=\Ver{\Q_1}\bigsqcup\cdots\bigsqcup\Ver{\Q_k}$ is
a disjoint union, for every $a\in\Arr{\Q}$ with $a',a''\in\Ver{\Q_i}$ for some $i$ we
have $a\in\Arr{\Q_i}$, and $k$ is the minimal number satisfying the given conditions.
%
%
%
Obviously, $\Inv(\Q,\n)$ is the tensor product of
$\Inv(\Q_1,\n_1),\ldots,\Inv(\Q_s,\n_s)$ for some dimension vectors $\n_1,\ldots,\n_s$ of
$\Q_1,\ldots,\Q_s$, respectively, satisfying $\n=\n_1\oplus\cdots\oplus\n_s$. Therefore,
it is sufficient to consider only strongly connected quivers.

Given a one-vertex quiver with $d$ loops ($d>2$), there are two possibilities for $D=D(\Q,(2))$:
\begin{enumerate}
\item[1)] if $\Char(K)=2$, then $D=d$;

\item[2)] if $\Char(K)\neq2$, then $D=3$.
\end{enumerate}
See Section~\ref{section_matrix_invariants} for the references. Note that in the first case $D$ depends linearly on $d$ and in the second case $D$ is a constant that does not depend on $d$. We show that the same statement is valid for an arbitrary quiver. Denote by $\QuiverNull(n,d,m)$ the set of all strongly connected quivers $\Q$ with
$\#\Ver{\Q}=n$, $\#\Arr{\Q}=d$, and $m(\Q)=m$. Our main result is the
following theorem.

\begin{theo}\label{theo_main_invariants}
Let $\Q\in\QuiverNull(n,d,m)$, where $d\geq2$. 
\begin{enumerate}
\item[1)] If $\Char(K)=2$, then $D(\Q,(2,\ldots,2))\leq md$. Moreover, if $n,m$ are fixed, then 
$$\frac{\max\{D(\Q,(2,\ldots,2))\,|\,\Q\in \QuiverNull(n,d,m)\}}{md}\to 1 \text{ as }d\to\infty.$$

\item[2)] If $\Char(K)\neq2$, then 
$D(\Q,(2,\ldots,2))\leq 3n$. Moreover, if $d$ is sufficiently large with respect to $n,m$, then $$\max\{D(\Q,(2,\ldots,2))\,|\,\Q\in \QuiverNull(n,d,m)\}$$%
is equal to the given bound.
\end{enumerate}
\end{theo} 

As an immediate corollary of Donkin's Theorem on the generators of $\Inv(\Q,\n)$ we obtain that the upper bounds on degrees from Theorem~\ref{theo_main_invariants} remain valid for the algebra of invariants $\Inv(\Q,(\de_1,\ldots,\de_n))$ with $\de_1,\ldots,\de_n\leq 2$.

The proof of Theorem~\ref{theo_main_invariants} is based on the reduction to the problem of finding out maximal paths
satisfying certain condition (see Lemma~\ref{lemma_reduction}). In informal way, the
last problem can be stated as ``Guide's Problem"{} (see below).

For a quiver $\Q$ introduce an equivalence $\equiv$ on the set of all closed paths
extended with an additional symbol $0$. For any paths $a,b$ such that $ab$ is a closed path and any incident closed paths $a_1,a_2,\ldots$ we
define
\begin{enumerate}
\item[1.] $ab\equiv ba$;

\item[2.] $a_{\si(1)}\cdots a_{\si(t)}\equiv (-1)^{\si}a_1\cdots a_t$, where $t\geq2$ and $\si\in
\Symm_t$;

\item[3.] $a_1^2a_2\equiv0$;

\item[4.] if $\Char(K)=2$, then $a_1^2\equiv0$; if $\Char(K)\neq2$, then
$a_1a_2a_3a_4\equiv0$.
\end{enumerate}

\begin{lemma}\label{lemma_reduction}
Let $a=a_1\cdots a_s$ be a closed path in $\Q$, where $a_1,\ldots,a_s\in\Arr{\Q}$. Then
$\tr(X_{a_s}\cdots X_{a_1})\in \Inv(\Q,(2,2,\ldots,2))$ is decomposable if and only if
$a\equiv0$.
\end{lemma}
\bigskip


Denote by $M(\Q)$ the maximal degree of a closed path $a$ in $\Q$ satisfying
$a\not\equiv0$. Lemma~\ref{lemma_reduction} shows that the case $\Char(K)=2$ is
essentially different from the case $\Char(K)\neq2$. The longest part of the paper is
dedicated to the case $\Char(K)=2$ and in this case $M(\Q)$ is equal to the length of a
route that provides a solution for the following problem.

\begin{max_path_problem}
{\it %
A guide shows a city to a tourist. They ride by car along streets of the city. All
streets are assumed to be one-way. (Two-way streets can be considered as two
different streets.) At the end of the tour they should come back to their starting point.
At the beginning the guide shows a plan of their route to the tourist. If the route goes
through a crossroad several times (i.e., this crossroad divides the route into parts
$ab_1\ldots b_kc$, where $k\geq1$ and $b_1,\ldots,b_s$ are cycles that start and
terminate at the given crossroad), then the tourist can choose order of passing these
cycles (i.e., the tourist can turn the route into $ab_{\pi(1)}\ldots b_{\pi(k)}c$ for any
permutation $\pi\in\Symm_k$). The route is called bad if it contains two consecutive
cycles that coincide.



Guide's payment depends on the length of their route and so his task is to find out the
longest route such that the tourist can not turn this route into a bad one.
}%
\end{max_path_problem}

In Section~\ref{section_matrix_invariants} we consider some results on generating systems
for invariants of a quiver with one vertex and several arrows. In Section~\ref{section_PR} we formulate Zubkov's Theorem, which we apply in Section~\ref{section_invariants} to prove
Lemma~\ref{lemma_reduction}. Section~\ref{section_notations} contains definitions of notions that are used in Sections~\ref{section_basic_equiv} and~\ref{section_example}.  If $\Char(K)\neq2$, then the upper bound on $M(\Q)$ is
calculated in Lemma~\ref{lemma_char_0}; otherwise, we establish the upper bound on
$M(\Q)$ in Corollary~\ref{cor_new}. In Lemma~\ref{lemma_example}
we estimate a precision of the given upper bounds. Taking into account
Lemma~\ref{lemma_reduction} and Remark~\ref{remark_det} together with the fact that
$\Inv(\Q,(2,2,\ldots,2))$ is generated by indecomposable invariants, we complete the
proof of Theorem~\ref{theo_main_invariants}. An example of indecomposable invariants is given in Example~\ref{ex_1}. 

\begin{remark}
In the next paper we will consider $n,d,m$ satisfying $\Q(n,d,m)\neq\emptyset$ and define the upper bound $M(n,d,m)$ such that in case $\Char(K)=2$ we have 
\begin{enumerate}
\item[$\bullet$] $D(\Q,(2,\ldots,2))\leq M(n,d,m)$ for all $\Q\in\Q(n,d,m)$; 

\item[$\bullet$] there is a $\Q\in\Q(n,d,m)$ such that $M(n,d,m)-m\leq D(\Q,(2,\ldots,2))$.
\end{enumerate}
\end{remark}

\section{Matrix invariants}\label{section_matrix_invariants}
Suppose $\Q$ is a quiver with one vertex and $d$ arrows. Then $\Inv(\Q,(n))$ is called
the {\it matrix} invariant algebra and we denote it by $R_{n,d}$. In this section we
discuss some known results on generating systems for $R_{n,d}$.

Relying on the theory of modules with good filtrations
(see~\cite{Donkin_1985},~\cite{Donkin_1993}), Donkin~\cite{Donkin_1992a} proved that
$K$-algebra $R_{n,d}\subset K[x_{ij}(r)\,|\,1\leq i,j\leq n,\,1\leq r\leq d]$ is
generated by $\si_k(X_{r_1}\cdots X_{r_s})$ for $1\leq k\leq n$ and $1\leq
r_1,\ldots,r_s\leq d$, where $X_r=(x_{ij}(r))_{1\leq i,j\leq n}$ is the $n\times n$
matrix. For zero characteristic case, generators were found earlier by Sibirskii
in~\cite{Sibirskii_1968} and Procesi in~\cite{Procesi_1976}.

A m.h.s.g.~for $R_{2,d}$ was found by Sibirskii in~\cite{Sibirskii_1968} when
$\Char(K)=0$, by Procesi in~\cite{Procesi_1984} when $\Char(K)$ is odd, and by Domokos,
Kuzmin, and Zubkov in~\cite{DKZ_2002} when $\Char(K)=2$. A m.h.s.g.~for $R_{3,d}$ was
found by the author in~\cite{Lopatin_Comm1},~\cite{Lopatin_Comm2}. Moreover, for $n=3$
and $d=2$ relations between elements of some m.h.s.g.~were explicitly described by
Nakamoto in~\cite{Nakamoto_2002} and by Aslaksen, Drensky, and Sadikova
in~\cite{ADS_2006} (see also Teranishi~\cite{Teranishi_1986}). A m.h.s.g.~for $R_{4,2}$
was described by Drensky and Sadikova in~\cite{Drensky_Sadikova_4x4} when $\Char(K)=0$.

Regarding an arbitrary $n$, an upper bound on indecomposable invariants of $R_{n,d}$ was
given by Domokos in~\cite{Domokos_gen_2002} in terms of the nilpotency degree $N(n,d)$ of
a (non-unitary) relatively free $d$-generated algebra with the identity $x^n=0$. If
$\Char(K)=0$ or $\Char(K)>n$, then $N(n,d)\leq 2^n-1$ by the Nagata--Higman Theorem
(see~\cite{Higman_1956}). Moreover, if $\Char(K)=0$, then $N(n,d)\leq n^2$ by Razmyslov
(see~\cite{Razmyslov_1974}) and $R_{n,d}$ is generated by elements of degree less or
equal to $N(n,d)$. The situation changes drastically when $0<\Char(k)\leq n$, namely,
$R_{n,d}$ is not generated by its elements of degree less than $d$. So if $n$ is fixed
and $d$ tends to infinity, then the maximal degree of indecomposable invariant as well as
$N(n,d)$ tends to infinity (see~\cite{DKZ_2002}). Observe that for an arbitrary
$\Char(K)$, there exists an upper bound on $N(n,d)$ by Klein (see~\cite{Klein_2000}):
$N(n,d)<(1/6)n^6d^n$. For more detailed introduction to finite generating systems for
$R_{n,d}$ see overviews~\cite{Formanek_1987} and~\cite{Formanek_1991} by Formanek. For
recent developments in characteristic zero see~\cite{Drensky_survey_2007} and in positive
characteristic see~\cite{DKZ_2002}.

\section{Zubkov's Theorem}\label{section_PR}

In what follows $\Q$ is a strongly connected quiver and $\n$ is its dimension vector. The
aim of this section is to formulate Zubkov's Theorem
(see~\cite{Zubkov_Fund_Math_2001}) that describes relations for the algebra of invariants
$\Inv(\Q,\n)$.

Denote by $S$ the free semigroup generated by letters $\{a_1,a_2,\ldots\}$. Words
$b=a_{i_1}\cdots a_{i_t}$ and $c=a_{j_1}\cdots a_{j_t}$ are called equivalent, if there
exists a cyclic permutation $\pi\in \Symm_t$ such that $i_k=j_{\pi(k)}$ for $1\leq k\leq
t$. The {\it cycle} (in letters $a_1,a_2,\ldots$) is the equivalence class of some word.
The cycle is {\it primitive}, if it is not equal to a power of a shorter cycle.

Let us recall some formulas. In this section $A,A_1,\ldots,A_s$ stand for $n\times n$
matrices and $n>1$. For $1\leq k\leq n$ Amitsur's formula states~\cite{Amitsur_1980}:
\begin{eq}\label{eq_Amitsur}
\si_k(A_1+\cdots+A_s)=\sum (-1)^{k-(j_1+\cdots+j_t)} \si_{j_1}(c_1)\cdots\si_{j_t}(c_t),
\end{eq}
where the sum ranges over all pairwise different primitive cycles $c_1,\ldots,c_t$ in
letters $A_1,\ldots,A_s$ and positive integers $j_1,\ldots,j_t$ with
$\sum_{i=1}^{t}j_i\deg(c_i)=k$. As an example,
$$
\si_2(A_1+A_2)=\si_2(A_1)+\si_2(A_2)+\si_1(A_1)\si_1(A_2)-\si_1(A_1A_2).
$$%
Denote the right hand side of~\Ref{eq_Amitsur} by $F_k(A_1,\ldots,A_s)$. Let $1\leq k\leq
n$ and $\al,\al_1,\ldots,\al_s\in K$. Using the formula
\begin{eq}
\si_k(\al A)=\al^k\si_k(A),
\end{eq}
we obtain
$$F_k(\al_1A_1,\ldots,\al_sA_s)=\sum \al^{\un{\de}}F_{\un{\de}}(A_1,\ldots,A_s),$$
where the sum ranges over all $\un{\de}=(\de_1,\ldots,\de_s)\in\NN^s$ with
$\de_1+\cdots+\de_s=k$, $\al^{\un{\de}}=\al_1^{\de_1}\cdots a_s^{\de_s}$, and
$F_{\un{\de}}(A_1,\ldots,A_s)$ is a polynomial in $\si_t(A_{i_1}\cdots A_{i_j})$. The
polynomial $F_{\un{\de}}(A_1,\ldots,A_s)$ is called a {\it partial linearization} of
$F_k(A_1,\ldots,A_s)$.

For $1\leq k\leq n$ and $r\geq2$ we have the following well-known formulas:
\begin{eq}\label{eq_A}
\si_k(A_1A_2)=\si_k(A_2A_1),
\end{eq}
\begin{eq}\label{eq_D}
\si_k(A^r)=\sum\limits_{i_1,\ldots,i_{kr}\geq0}\beta^{(k,r)}_{i_1,\ldots,i_{kr}}
    \si_1(A)^{i_1}\cdots\si_{kr}(A)^{i_{kr}},%
\end{eq}
\noindent where we assume that $\si_i(A)=0$ for $i>n$. Denote the right hand side
of~\Ref{eq_D} by $G_{k,r}(A)$. In~\Ref{eq_D} coefficients
$\beta^{(k,r)}_{i_1,\ldots,i_{kr}} \in \ZZ$ do not depend on $A$ and $n$. If we take
$A=\diag(a_1,\ldots,a_n)$ is a diagonal matrix, then $\si_k(A^r)$ ($\si_i(A)$ for $1\leq
i\leq n$, respectively) is a symmetric polynomial (the $i$-th elementary symmetric
polynomial, respectively) in $a_1,\ldots,a_n$ and the coefficients
$\beta^{(k,r)}_{i_1,\ldots,i_{kr}}$ with $i_1+2 \cdot i_2+\cdots +kr\cdot i_{kr}\leq n$
can easily be found. As an example,
\begin{eq}\label{eq_tr_a2}
\si_1(A^2)=\si_1(A)^2-2\si_2(A).
\end{eq}

Consider the commutative algebra $A(\Q)$, freely generated by ``symbolic"{} elements
$\si_k(h)$, where $k\geq1$ and $h$ is a closed path in $\Q$. An expression $p=q$ with
different $p,q\in A(\Q)$ we interpret as the element $p-q$ of $A(\Q)$.

By Zubkov's Theorem, the algebra of invariants $\Inv(\Q,\n)$ is isomorphic
to $A(\Q)/T(\Q,\n)$, where the ideal $T(\Q,\n)$ is generated by
\begin{enumerate}
\item[$(A)$] $\si_k(a_1a_2)=\si_k(a_2a_1)$, where $k\geq1$ and $a_1,a_2$ are such paths
in $\Q$ that $a_1a_2$ is a closed path;

\item[$(B)$] $\si_k(a^r)=G_{k,r}(a)$, where $k\geq1$, $r\geq2$ and $a$ is a closed path
in $\Q$;

\item[$(C)$] $F_{\un{\de}}(a_1,\ldots,a_s)=0$, where $a_1,\ldots,a_s$ are incident
closed paths in $\Q$ for $s\geq1$, $\un{\de}=(\de_1,\ldots,\de_s)\in\NN^s$, and
$\de_1+\cdots+\de_s>\n_{a_1''}$; in particular, $\si_k(a)=0$ for any closed path $a$ in $\Q$ and $k>\n_{a''}$.
\end{enumerate}
The isomorphism is given by
$$ \si_k(X_{a_s}\cdots X_{a_1})\to\si_k(a_1\cdots a_s)$$%
for $1\leq k\leq \n_{a_1''}$ and such arrows $a_1,\ldots,a_s\in\Arr{\Q}$ that $a_1\cdots
a_s$ is a closed path in $\Q$. Elements of $T(\Q,\n)$ are called {\it relations} for
$\Inv(\Q,\n)$.

\section{Relations between indecomposable invariants}\label{section_invariants}

In this section $\Q$ is a strongly connected quiver. For short, denote
$\Inv(\Q,(2,2,\ldots,2))$ by $\Inv(\Q)$ and $T(\Q,(2,2,\ldots,2))$ by $T(\Q)$. If $a$ and $b$ are equal elements of
$\ov{\Inv(\Q)}={\Inv(\Q)}/{(\Inv(\Q)^{+})^2}$, then we write $a\equiv b$.

Consider the algebra $A(\Q)$ that was defined in Section~\ref{section_PR} and denote:
$$\tr(a)=\si_1(a)\text{ and }\det(a)=\si_2(a)$$
for a closed path $a$ in $\Q$. The algebra $\ov{\Inv(\Q)}$ is isomorphic to
$\ov{A(\Q)}/L(\Q)$ for $\ov{A(\Q)}={A(\Q)}/{(A(\Q)^{+})^2}$ and some ideal $L(\Q)\lhd
\ov{A(\Q)}$. As above, an element $q\in L(\Q)$ is called a relation for $\ov{\Inv(\Q)}$
and we write $q\equiv0$. We say that $q\equiv0$ follows from
$q_1\equiv0,\ldots,q_s\equiv0$, if $q$ is a linear combination of $q_1,\ldots,q_s$.

\begin{lemma}\label{lemma_relations}
The ideal of relations for $\ov{\Inv(\Q)}$ is equal to the $K$-span of the relations:
\begin{enumerate}
\item[(a)] $\tr(a_1a_2)\equiv \tr(a_2a_1)$, $\det(a_1a_2)\equiv \det(a_2a_1)$, where
$a_1$ and $a_2$ are such paths in $Q$ that $a_1a_2$ is a closed path;

\item[(b)] $\tr(a^2)\equiv-2\det(a)$;

\item[(c)] $\tr(a^2b)\equiv0$;

\item[(d)] $\tr(bac)\equiv-\tr(abc)$, $\tr(acb)\equiv-\tr(abc)$;

\item[(e)] $\det(ab)\equiv0$;

\item[(f)] $\si_k(a)\equiv0$ for $k>2$;
\end{enumerate}
where $a$, $b$, and $c$ are incident closed paths in $\Q$.
\end{lemma}
\begin{proof} Denote by $(\ov{A})$ relations for $\ov{\Inv(\Q)}$ obtained by factorization
of $(A)$ modulo the ideal $(A(\Q)^{+})^2$. Similarly define $(\ov{B})$ and $(\ov{C})$:
\begin{enumerate}
\item[$(\ov{B})$] $\si_k(a^r)\equiv \al_{k,r}\si_{kr}(a)$, where $k\geq1$, $r\geq2$, $a$
is a closed path in $\Q$, and $\al_{k,r}\in\ZZ$ do not depend on $a$;

\item[$(\ov{C})$] $\ov{F}_{\un{\de}}(a_1,\ldots,a_s)\equiv0$, where $a_1,\ldots,a_s$ are
incident closed paths in $\Q$, $\un{\de}=(\de_1,\ldots,\de_s)\in\NN^s$, and
$\de_1+\cdots+\de_s>2$.
\end{enumerate}
Let $a_1,\ldots,a_s$ be such arrows of $\Q$ that $a=a_1\cdots a_s$ is a closed path in $\Q$ and let $k=1,2$. Obviously, $\si_k(X_{a_s}\cdots X_{a_1})$ is decomposable if and only if $\si_k(a)\in T(\Q)+(A(\Q)^{+})^2$. Consider the grading on the polynomial ring $A(\Q)$ such that the generators have degree one. Then $\si_k(a)$ is homogeneous of degree one. Denote by $T(\Q)_1\subset A(\Q)$ the vector space consisting of the linear components of the elements of $T(\Q)$. We conclude that $\si_k(X_{a_s}\cdots X_{a_1})$ is decomposable if and only if $\si_k(a)\in T(\Q)_1$. The same reasoning shows that the image of $q\in A(\Q)$ in $\Inv(\Q)\simeq A(\Q)/T(\Q)$ under the canonical homomorphism is decomposable if and only if $q\in T(\Q)_1$. Thus, relations for $\ov{\Inv(\Q)}$ are linear combinations of relations $(\ov{A})$, $(\ov{B})$, and $(\ov{C})$ (see Section~\ref{section_PR}).

\smallskip
\textbf{1.} Let us prove that $(a)$--$(f)$ are relations for $\ov{\Inv(\Q)}$. First of
all notice that $(a)$ follows from $(\ov{A})$ and $(e)$ follows from $(\ov{C})$.
By~\Ref{eq_tr_a2}, the relation~$(\ov{B})$ with $k=1$ and $r=2$ coincides with~$(b)$.

If $\un{\de}=(2,1)$, then $(\ov{C})$ coincides with~$(c)$. If $\un{\de}=(2,2)$, then
$(\ov{C})$ is $\det(a_1a_2)-\tr(a_1^2a_2^2)\equiv0$. Since~$(c)$ implies
$\tr(a_1^2a_2^2)\equiv0$, $(d)$ is a relation. If $\un{\de}=(1,1,1)$, then $(\ov{C})$ is
$\tr(a_1a_2a_3)+\tr(a_2a_1a_3)\equiv0$. Clearly,~$(d)$ follows from this relation
and~$(a)$.

\smallskip
\textbf{2.} Let us prove that relations $(\ov{A})$, $(\ov{B})$, and $(\ov{C})$ follow
from $(a)$--$(f)$. Note that $(a)$, $(c)$, and $(d)$ imply that for all incident closed
paths $a$, $b$, and $c$ we have
\begin{eq}\label{eq_tr_abac}
\tr(abac)\equiv0.
\end{eq}
\indent Obviously, $(\ov{A})$ follows from~$(a)$ and~$(e)$.

Consider the relation $(\ov{B})$ with $k=1$. If $r=2$, then $(\ov{B})$ coincides
with~$(b)$. If $r\geq 3$, then $\si_r(a)\equiv0$ follows from~$(f)$ and $\tr(a^r)\equiv0$
follows from~$(c)$; hence, $(\ov{B})$ follows from~$(c)$ and~$(f)$.

If $k=2$ and $r\geq2$, then $\det(a^r)\equiv0$ follows from~$(e)$ and
$\si_{2r}(a)\equiv0$ follows from~$(f)$, so $(\ov{B})$ follows from~$(e)$ and~$(f)$. If
$k\geq3$, then $(\ov{B})$ follows from~$(f)$.

Now we consider the relation $(\ov{C})$. If $s=1$, then $(\ov{C})$ coincides
with~$(f)$. If $\de_i\geq2$ for some $i$, then $(\ov{C})$ follows from~$(e)$, $(f)$,
and~\Ref{eq_tr_abac}. If $\un{\de}=(1,\ldots,1)$, then~$(\ov{C})$ is
$$
\sum_{\pi\in\Symm_{s-1}}\tr(a_{\pi(1)}\cdots a_{\pi(s-1)}a_{s})\equiv0.
$$
The last relation follows from~$(a)$ and~$(d)$.
\end{proof}

Now we can prove Lemma~\ref{lemma_reduction} (see Section~\ref{section_intro}).

\bigskip
\begin{proof_lemma_reduction}
We have the following relations for $\Inv(\Q)$:
\begin{enumerate}
\item[$(g)$] $\tr(a_{\si(1)}\cdots a_{\si(t)})\equiv \sign(\si)\tr(a_1\cdots a_t)$ for
$t\geq1$ and $\si\in \Symm_t$;

\item[$(h)$] if $\Char(K)\neq2$, then $\tr(a_1\cdots a_4)\equiv0$;
\end{enumerate}
where $a_1,a_2,\ldots$ are incident closed paths in $\Q$. To prove this notice that $(g)$
follows from~$(a)$ and~$(d)$ (see Lemma~\ref{lemma_relations}). Consecutively using
relations~$(d)$ and~$(g)$ we obtain:
$$
\tr(a_1a_2a_3a_4)\equiv - \tr(a_3a_1a_2a_4)\equiv -\tr(a_1a_2a_3a_4).
$$
Hence $(h)$ is a relation for $\Inv(\Q)$.

Suppose $\tr(a)\equiv0$. By Lemma~\ref{lemma_relations}, the relation $\tr(a)\equiv0$
follows from~$(a)$--$(f)$. Hence this relation follows from $(a)$, $(b)$, $(g)$ or from $(a)$, $(c)$, $(g)$ when $\Char(K)=2$; and it follows from $(a)$, $(c)$, $(g)$, $(h)$ when $\Char(K)\neq2$. In all cases we have $a\equiv0$.

If $a\equiv0$, then Lemma~\ref{lemma_relations} together with~$(g)$ and~$(h)$ gives
$\tr(a)\equiv0$.
\end{proof_lemma_reduction}

%
%
%
%

\begin{remark}\label{remark_det}
Let $a=a_1\cdots a_s$ be a closed path in $\Q$, where $a_1,\ldots,a_s\in\Arr{\Q}$. If
$q=\det(X_{a_s}\cdots X_{a_1})\in \Inv(\Q)$ is indecomposable, then $\deg(q)\leq
2m$.
\end{remark} 
\begin{proof}
If $q$ is indecomposable, then $a$ is a primitive closed path and $\deg(a)\leq m$.
\end{proof}

\section{Some notations and auxiliary results}\label{section_notations}
Suppose $a=a_1\cdots a_s$ is a path in a quiver $\Q$ and $a_1,\ldots,a_s\in\Arr{\Q}$.
Given $v\in\Ver{\Q}$ and $b\in\Arr{\Q}$, let $\deg_b(a)=\#\{i\,|\,a_i=b,\,1\leq i\leq
s\}$ be the degree of $a$ in the arrow $b$ and $\deg_v(a)=\#\{i\,|\,a_i'=v,\,1\leq i\leq
s\}+\rho$ be
the degree of $a$ in the vertex $v$, where %
$$
\rho=\left\{
\begin{array}{rl}
1,& \text{if } a_1''=v\text{ and } a_s'\neq v\\
0,&\text{otherwise}\\
\end{array}
\right..
$$
As an example, if $a_1'',a_1',\ldots,a_s'$ are pairwise different, then
$\deg_{a_1''}(a)=\deg_{a_1'}(a)=\cdots=\deg_{a_s'}(a)=1$. If $a$ is known to be a closed
path, then $\deg_v(a)=\#\{i\,|\,a_i'=v,\,1\leq i\leq s\}$ coincides with the definition
given in Section~\ref{section_intro}. Sometimes it is convenient to consider only those
vertices of $a$ that are not equal to $a'$ and $a''$; for this purpose introduce
$\degII{v}{a}=\#\{i\,|\,a_i'=v,\,1\leq i\leq s-1\}$. Note that $\deg_v(a)=\degII{v}{a}$
for $v\not\in\{a',a''\}$.

The {\it multidegree} of a path $a$ in $\Q$ is $\un{\de}=(\de_b)_{b\in\Arr{\Q}}$, where
$\de_b=\deg_b(a)$, and we denote it by $\mdeg(a)$. 

Let $x_1,\ldots,x_s$ be all arrows in $\Q$ from $u$ to $v$, where $u,v\in\Ver{\Q}$. Then
denote by $\ch{x}$ any arrow from $x_1,\ldots,x_s$, by $\{\ch{x}\}$ the set
$\{x_1,\ldots,x_s\}$, and say that $\ch{x}$ is an arrow from $u$ to $v$. Schematically,
we depict arrows $x_1,\ldots,x_s$ as
$$
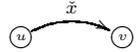
\vcenter{
\xymatrix@C=1cm@R=1cm{ %
\vtx{u}\ar@/^/@{->}[r]^{\ch{x}} &\vtx{v}\\
}}.
$$
For a path $a$ in $\Q$ denote $\deg_{\ch{x}}(a)=\sum_{i=1}^s\deg_{x_i}(a)$. As an
example, an expression $\ch{x}a_1\cdots\ch{x}a_k$ stands for a path $x_{i_1}a_1\cdots
x_{i_k}a_k$ for some $1\leq i_j\leq s$ ($1\leq j\leq k$). Similarly, if $x_1,\ldots,x_s$
are loops in $v\in\Ver{\Q}$, then $\ch{x}^k$ stands for a closed path $x_{i_1}\cdots
x_{i_k}$ for some $i_1,\ldots,i_k$.

Suppose $a=a_1\cdots a_s$ is a primitive closed path in $\Q$ and
$a_1,\ldots,a_s\in\Arr{\Q}$. The path $a$ is called {\it $\un{\de}$-single} if
$\de_{a_i}\geq1$ for all $i$ and $\de_{a_i}=1$ for some $i$. The path $a$ is called {\it
$\un{\de}$-double} if $\de_{a_i}\geq2$ for all $i$.

For a $v\in\Ver{\Q}$ denote by $1_v$ the {\it empty path} in the vertex $v$. Given a path
$b$ with $b'=v$, we write $b1_v=b$ and for a path $b$ with $b''=v$ we write $1_vb=b$. By
definition, $\deg(1_v)=0$, $\deg_a(1_v)=0$, and
$$\deg_w(1_v)=\left\{
\begin{array}{rl}
1,& \text{if }  w=v\\
0,& \text{if } w\neq v\\
\end{array}
\right.
$$
for all $w\in\Ver{\Q}$ and $a\in\Arr{\Q}$. Denote by $\path(\Q)$ the set of all paths and
empty paths in $\Q$. If we consider a path, then we assume that it is non-empty unless
otherwise stated; if we write $a\in\path(\Q)$, then we assume that a path $a$ can be
empty.

For closed paths $a,b$ we write $a\sim b$ if $a=c_1c_2$ and $b=c_2c_1$ for some
$c_1,c_2\in\path(\Q)$. A path $b$ is called a {\it subpath} in a path $a$, if the path
$a$ is closed and $a\sim bc$, or $a$ is not closed and $a=c_1bc_2$, where
$c,c_1,c_2\in\path(\Q)$.

Suppose  $V\subset\Ver{\Q}$ is a subset and a path $h$ in $\Q$ satisfies $h',h''\in V$.
We say that a quiver $\QG$ is the {\it $h$-restriction} of $\Q$ to $V$ if $\Ver{\QG}=V$
and $\Arr{\QG}=\{\tilde{a}\}$, where $a$ ranges over such subpaths of $h$ that $a',a''\in
V$ and $\degII{v}{a}=0$ for all $v\in V$. By definition, $\tilde{a}'=a'$ and
$\tilde{a}''=a''$. There is a unique path in $\QG$ that corresponds to $h$ and each path
in $\QG$ corresponds to some path in $\Q$.

\begin{example} Let $\Q$ be the quiver
$$ \vcenter
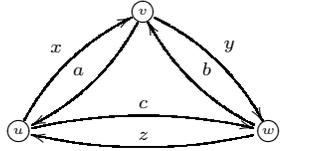
{
\xymatrix@C=1.3cm@R=1.3cm{ %
&\vtx{v} %
\ar@/^/@{->}[ld]_{a} \ar@/_/@{<-}[ld]_{x}
\ar@/^/@{->}[rd]^{y} \ar@/_/@{<-}[rd]^{b}& \\
\vtx{u} \ar@/^/@{->}[rr]^{c} \ar@/_/@{<-}[rr]^{z} &&\vtx{w}  }}%
\qquad,
$$
$h$ be a path in $\Q$ with $h',h''\in\{u,v\}$, and $\QG$ be the $h$-restriction of $\Q$
to the vertices $u$ and $v$. Then $\QG$ is a subquiver of the quiver
$$
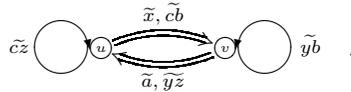
\loopR{0}{0}{\widetilde{cz}}%
\xymatrix@C=1.3cm@R=1.3cm{ %
\vtx{u}\ar@2@/^/@{->}[r]^{\widetilde{x},\,\widetilde{cb}}\ar@2@/_/@{<-}[r]_{\widetilde{a},\,\widetilde{yz}} &\vtx{v} \\
}%
\loopL{0}{0}{\,\widetilde{yb}}\qquad\,.
$$
\end{example}

Dealing with equivalences we use the following conventions. If we write $a\equiv b$,
then we assume that $a$ and $b$ are closed paths in $\Q$. If we write $ab$ for paths $a$
and $b$, then we assume that $a'=b''$. To explain how we apply formulas to prove some
equivalence $a\equiv b$ we split the word $a$ into parts using dots. As an example, see
the proof of part~1 of Lemma~\ref{lemma_star}.

The next lemma is well known.

\begin{lemma}\label{lemma_mdeg} 
Suppose $Q$ is a strongly connected quiver and $\un{\de}\in\NN^{\#\Arr{\Q}}$. Then the
following conditions are equivalent:
\begin{enumerate}
\item[a)] There is a closed path $h$ in $\Q$ such that $\mdeg(h)=\un{\de}$ and
$\Arr{h}=\Arr{\Q}$; in particular, $\Ver{h}=\Ver{\Q}$.

\item[b)] We have $\de_a\geq1$ for all $a\in\Arr{\Q}$ and $\sum_{a'=v}\de_a =
\sum_{a''=v}\de_a$ for all $v\in\Ver{\Q}$, where the sums range over all $a\in\Arr{\Q}$
satisfying the given conditions.
\end{enumerate}
\end{lemma}

\section{Upper bounds}\label{section_basic_equiv}

Let $\Q$ be a quiver. We start with the case of $\Char(K)\neq2$.

\begin{lemma}\label{lemma_char_0}
Suppose $\Char(K)\neq2$. If $\Q$ is a quiver with $n$ vertices and $h$ is a closed path in $\Q$ with $h\not\equiv0$, then $\deg(h)\leq 3n$.
\end{lemma}
\begin{proof}
If the claim of the lemma is wrong, then there is a vertex $v\in\Ver{\Q}$ such that $\deg_v(h)\geq 4$. Then $h\equiv h_1\cdots h_4$ for some closed paths $h_1,\ldots,h_4$ in $v$. Thus $h\equiv0$ by the definition of the equivalence $\equiv$; a contradiction.
\end{proof}

In what follows we assume $\Char(K)=2$ unless otherwise stated.  We will use the following remark without references to it.

\begin{remark}\label{remark_no_change}
Suppose $f,h$ are closed paths in $\Q$ and $b$ is a subpath of $f$. Let the equivalence
$f\equiv h$ follows from the formulas of the form $a_{\si(1)}\cdots a_{\si(t)}\equiv
a_1\cdots a_t$, where $a_1,\ldots,a_t$ are closed paths in $v\in\Ver{\Q}$ satisfying
$\degII{v}{b}=0$, $t\geq2$, and $\si\in \Symm_t$. Then $b$ is also a subpath of $h$.
\end{remark}

\begin{lemma}\label{lemma_L0}
Let $h$ be a closed path in $\Q$ and $\{\ch{p}\}$ be loops of $\Q$ in some
$v\in\Ver{\Q}$. Then $h\equiv \ch{p}^kb$, where $k\geq0$, $b\in\path(\Q)$, and
$\deg_{\ch{p}}(b)=0$.

Moreover, suppose $a\in\Arr{h}$ and $a'\neq a''$. If $a'=v$, then $h\equiv a
\ch{p}^kb_0$; if $a''=v$, then $h\equiv \ch{p}^kab_0$, where, as above,
$\deg_{\ch{p}}(b_0)=0$.
\end{lemma}
\begin{proof}
Denote $k=\deg_{\ch{p}}(h)$. If $k\leq1$, then the statement of the lemma is trivial.
Otherwise $h\sim \ch{p}g_1\cdots \ch{p}g_k\equiv \ch{p}^kg_1\cdots g_k$ for some
$g_1,\ldots,g_k\in\path(\Q)$. The proof of the second part of the lemma is similar.
\end{proof}

\begin{lemma}\label{lemma_star}
1. We have $x_1 a_1 x_2 a_2 x_3 a_3\equiv x_3 a_1 x_1 a_2 x_2 a_3$, where $x_1,x_2,x_3$
are paths from $u$ to $v$ and $a_1,a_2,a_3$ are paths from $v$ to $u$ for
$u,v\in\Ver{\Q}$.

2. Let $h=\ch{x} a_1\cdots \ch{x} a_s$ be a closed path.
\begin{enumerate}
\item[a)] If $\deg_{\ch{x}}(h)\geq 3$ and $\deg_{x_1}(h)\geq 1$, then %
$h\equiv x_1 a_1 \ch{x} a_2 \cdots \ch{x} a_s$.

\item[b)] If $\deg_{x_1}(h)\geq 2$, then %
$h\equiv x_1 a_1 x_1 a_2\ch{x} a_2 \cdots \ch{x} a_s$.

\item[c)] If $\deg_{x_1}(h)\geq 1$ and $\deg_{x_2}(h)\geq 1$, then %
$h\equiv x_1 a_1 x_2 a_2\ch{x} a_2 \cdots \ch{x} a_s$ or $h\equiv x_2 a_1 x_1 a_2\ch{x}
a_2 \cdots \ch{x} a_s$.
\end{enumerate}
\end{lemma}
\begin{proof}
\textbf{1. } Equivalences %
$x_1\cdot a_1 x_2\cdot a_2 x_3\cdot a_3\equiv %
x_1\cdot a_2 x_3\cdot a_1 x_2\cdot a_3 =%
x_1 a_2\cdot x_3 a_1 \cdot x_2 a_3 \equiv %
x_3 a_1 x_1 a_2 x_2 a_3$ give the required formula.

\textbf{2a) } There are three possibilities: %
$h\sim x_i c_1 x_1 c_2 x_j c_3$, %
$h\sim x_i c_1 x_j c_2 x_1 c_3$, and %
$h\sim x_1 c_1 x_i c_2 x_j c_3$ for some paths $c_1,c_2,c_3$ and numbers $i,j$. Applying
part~1, we obtain the claim. Similarly, we prove parts~2b) and~2c).
\end{proof}

\begin{lemma}\label{lemma_Q_with 2ver}
Let $h$ be a closed path in a quiver $\Q$
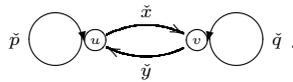
$$\loopR{0}{0}{\ch{p}\,}
\xymatrix@C=1cm@R=1cm{ %
\vtx{u}\ar@/^/@{->}[r]^{\ch{x}}\ar@/_/@{<-}[r]_{\ch{y}} &\vtx{v} \\
}\loopL{0}{0}{\,\ch{q}}\quad.
$$
\begin{enumerate}
\item[a)] If $\deg_{\ch{x}}(h)+\deg_{\ch{y}}(h)\geq1$, then $h\equiv
\ch{p}^i\,\ch{x}\,\ch{q}^j\,\ch{y} (\ch{x}\ch{y})^k$ for some $i,j,k\geq0$.

\item[b)] If $\deg_{x_1}(h)\geq2$, $\deg_{y_1}(h)\geq2$, and
$\deg_{\ch{x}}(h)+\deg_{\ch{y}}(h)>4$, then $h\equiv (x_1y_1)^2f$ for some path $f$.
\end{enumerate}
\end{lemma}
\begin{proof} \textbf{a)} Using Lemma~\ref{lemma_L0}, we have %
$h\equiv \ch{p}^i\,a\,\ch{q}^j\,b$ for some paths $a,b$ and $i,j\geq0$. Equalities
$a=\ch{x}\,(\ch{y}\ch{x})^k$ and $b=\ch{y}\,(\ch{x}\ch{y})^l$ for $k,l\geq0$ complete the
proof.

\smallskip
\textbf{b)}  Part~a) implies that $h\equiv \ch{p}^ix_l\ch{q}^jf$ for $f=\ch{y}
(\ch{x}\ch{y})^k$, where $i,j,k\geq0$ and $l\geq1$. If $\deg_{\ch{x}}(h)\geq3$, then,
taking into account part~1 of Lemma~\ref{lemma_star}, we can assume that
$\deg_{x_1}(f)\geq2$. We add a new arrow $x_0$ to $\Q$ and define $x_0'=v$, $x_0''=u$.
Using part~a) of the lemma together with part~2b) of Lemma~\ref{lemma_star}, we obtain
the required equivalence for the closed path $x_0f$ . Substituting $\ch{p}^ix_l\ch{q}^j$
for $x_0$, we prove the required equivalence for $h$. The case $\deg_{\ch{y}}(h)\geq3$ is
similar.
\end{proof}

\begin{lemma}\label{lemma_Q_with 3ver}
Suppose $\Q$ is
$$
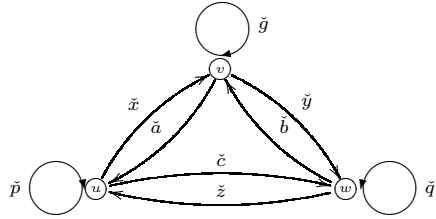
\loopD{-83}{-8}{\,\ch{g}} \loopR{0}{45}{\ch{p}\,}
\xymatrix@C=1.3cm@R=1.3cm{ %
&\vtx{v} %
\ar@/^/@{->}[ld]_{\ch{a}} \ar@/_/@{<-}[ld]_{\ch{x}}
\ar@/^/@{->}[rd]^{\ch{y}} \ar@/_/@{<-}[rd]^{\ch{b}}& \\
\vtx{u} \ar@/^/@{->}[rr]^{\ch{c}} \ar@/_/@{<-}[rr]^{\ch{z}} &&\vtx{w}  }
\loopL{-1}{45}{\,\ch{q}} %
$$
\smallskip
and $h$ is a closed path with $\deg_{x_1}(h)\geq2$ and $\deg_{y_1}(h)\geq2$.
\begin{enumerate}
\item[a)] If $\deg_{\ch{g}}(h)=0$, then $h\equiv x_1y_1 f_1 x_1y_1 f_2$ for some paths
$f_1$ and $f_2$.

\item[b)] If $\deg_{\ch{x}}(h)+\deg_{\ch{y}}(h)+\deg_{\ch{a}}(h)+\deg_{\ch{b}}(h)>4$,
then $h\equiv x_1y_1 f_1 x_1y_1 f_2$ for some paths $f_1$ and $f_2$.
\end{enumerate}
\end{lemma}
\begin{proof} \textbf{a)}
For short, in this proof we use one and the same symbol $l$ for non-negative integers
that can be different. As an example, $\ch{p}^l\ch{c}\,\ch{q}^l$ stands for
$\ch{p}^{l_1}\ch{c}\,\ch{q}^{l_2}$ for some $l_1,l_2\geq0$. Let $\QG$ be the
$h$-restriction of $\Q$ to the vertices $u$ and $v$:
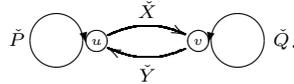
$$\loopR{0}{0}{\ch{P}}
\xymatrix@C=1cm@R=1cm{ %
\vtx{u}\ar@/^/@{->}[r]^{\ch{X}}\ar@/_/@{<-}[r]_{\ch{Y}} &\vtx{v} \\
}\loopL{0}{0}{\,\ch{Q}}\quad.
$$
We have the inclusions
$$\{\ch{X}\}\subset \{\ch{x},\ch{c}\ch{q}^l\ch{b}\},\,%
\{\ch{Y}\}\subset \{\ch{a},\ch{y}\ch{q}^l\ch{z}\},\, %
\{\ch{P}\}\subset \{\ch{p},\ch{c}\ch{q}^l\ch{z}\},\, %
\{\ch{Q}\}\subset \{\ch{y}\ch{q}^l\ch{b}\}.
$$
Consider $h$ as a path in $\QG$. Part~a) of Lemma~\ref{lemma_Q_with 2ver} together with
$\deg_{\ch{X}}(h)\geq \deg_{\ch{x}}(h)\geq2$ implies that %
$h\equiv \ch{P}^i\,\ch{X}\,\ch{Q}^j\,\ch{Y} (\ch{X}\ch{Y})^k$ for some $i,j\geq0$ and
$k\geq1$.
Moreover, applying part~2b) of Lemma~\ref{lemma_star} to $\ch{X}$, we obtain %
$h\equiv \ch{P}^i\,x_1\,\ch{Q}^j\,\ch{Y}\, x_1\ch{Y}(\ch{X}\ch{Y})^{k-1}$.

If $j=0$, then $\{\ch{Y}\}$ contains $Y_1=y_1\ch{q}^{l_1}\ch{z}$ and
$Y_2=y_1\ch{q}^{l_2}\ch{z}$ for some $l_1,l_2\geq0$. We apply part~2c) of
Lemma~\ref{lemma_star} to $\ch{Y}$ and obtain %
$h\equiv \ch{P}^i\,x_1\,Y_{s_1} x_1 Y_{s_2}(\ch{X}\ch{Y})^{k-1}$, where
$s_1,s_2\in\{1,2\}$ and $s_1\neq s_2$. The required equivalence is proven.

If $j\geq1$, then we rewrite $h$ as %
$$h\equiv \ch{P}^i\,x_1\,y_s\ch{q}^l\ch{b}\,R\,\ch{Y}\, x_1Y_1(\ch{X}\ch{Y})^{k-1}$$
for
$s\geq1$, $R=\ch{Q}^{j-1}$, and $Y_1\in\{\ch{Y}\}$.

If $\{\ch{Y}\}\subset \{\ch{a}\}$, then $j\geq2$. Applying part~2b) of
Lemma~\ref{lemma_star} to $\ch{y}$, we can assume that $s=1$ and
$R=y_1\ch{q}^l\ch{b}\,\ch{Q}^{j-2}$. Hence
$$
h\equiv \ch{P}^i\,x_1\cdot %
y_1\ch{q}^l\ch{b}\cdot %
y_1\ch{q}^l\ch{b}\,\ch{Q}^{j-2}\,\ch{Y}\,x_1\cdot %
\ch{Y}(\ch{X}\ch{Y})^{k-1}$$ %
$$\quad\,\equiv
\ch{P}^i\,x_1\cdot %
y_1\ch{q}^l\ch{b}\,\ch{Q}^{j-2}\,\ch{Y}\,x_1\cdot %
y_1\ch{q}^l\ch{b}\cdot %
\ch{Y}(\ch{X}\ch{Y})^{k-1}.$$ %

If $\{\ch{Y}\}$ is not a subset of $\{\ch{a}\}$, then $\{\ch{Y}\}$ contains
$y_r\ch{q}^l\ch{z}$ for some $r\geq1$. Applying part~2b) of Lemma~\ref{lemma_star} to
$\ch{y}$, we can assume that $s=r=1$. There are three cases:
\begin{enumerate}
\item[1.] If $Y_1=y_1\ch{q}^l\ch{z}$, then the claim is proven.

\item[2.] If $k\geq2$, then $\deg_{\ch{Y}}(h)=k+1\geq3$. Applying part~2a) of
Lemma~\ref{lemma_star} to $\ch{Y}$, we can assume that $Y_1=y_1\ch{q}^l\ch{z}$ and the
claim is proven.

\item[3.] If $k=1$ and $Y_1\neq y_1\ch{q}^l\ch{z}$, then
$$h\equiv
\ch{P}^i\,x_1\cdot %
y_1\ch{q}^l\ch{b}\,R\cdot%
y_1\ch{q}^l\ch{z}\, x_1\cdot %
Y_1%
\equiv \ch{P}^i\,x_1\cdot %
y_1\ch{q}^l\ch{z}\, x_1\cdot %
y_1\ch{q}^l\ch{b}\,R\cdot%
Y_1.$$
\end{enumerate}

Part~a) of the lemma is proven.

\smallskip
\textbf{b)} If $\deg_{\ch{a}}(h)\geq1$, then, taking into account Lemma~\ref{lemma_L0},
we have $h\equiv \ch{g}^k\ch{a}f$ for $k\geq0$ and a path $f$ with $\deg_{\ch{g}}(f)=0$.
We add a new arrow $a_0$ to $\Q$ and define $a_0'=u$, $a_0''=v$. Then the closed path
$a_0f$ satisfies the condition of part~a) of the lemma and the required equivalence is
valid for it. Substituting $\ch{g}^k\ch{a}$ for $a_0$, we prove the required equivalence
for $h$. The remaining cases $\deg_{\ch{b}}(h)\geq1$, $\deg_{\ch{x}}(h)\geq3$, and
$\deg_{\ch{y}}(h)\geq3$ can be treated analogously.
\end{proof}

\begin{remark} Note that the conditions from parts~a),~b)
can not be omitted. As an example, if $h=g_1y_1z_1x_1y_1z_1x_1$, then for any paths
$f_1,f_2$ we have $h\not\equiv x_1y_1f_1x_1y_1f_2$.
\end{remark}
\bigskip

Suppose a quiver $\Q$ contains a path $a=a_1\cdots a_s$, where
$a_1,\ldots,a_s\in\Arr{\Q}$ are pairwise different.  Let $h$ be a closed path in $\Q$ such that $\deg_{a_i}(h)\geq2$ for all $i$
and there is a $b\in\Arr{h}$ satisfying $b\neq a_i$ for all $i$.


\begin{lemma}\label{lemma_aaa}
Using the preceding notation we have $h\equiv a_1\cdots a_sf$ for some $f\in\path(\Q)$.
Moreover,
\begin{enumerate}
\item[a)] if $b'=a''_1$, then $h\equiv b\,a_1\cdots a_sf$ for some $f\in\path(\Q)$;

\item[b)] if $b''=a'_s$, then $h\equiv a_1\cdots a_s\,bf$ for some $f\in\path(\Q)$.
\end{enumerate}
\end{lemma}
\begin{proof}
Let us prove part~a). Since $\deg_{a_1}(h)\geq2$ and $b\in\Arr{h}$, we have $h\sim
gb\cdot q\cdot a_1 f$ for paths $g,q,f$, where $q'=q''=v_1$, $g''=f'=v_1$, and $q$ can be
empty.

If $q$ is empty, then $h\sim gb\,a_1 f$. If $q$ is non-empty, then $h\equiv gb\cdot a_1
f\cdot q$. Thus in both cases we have $h\equiv b\,a_1f_1$ for some path $f_1$. Continuing
this procedure we complete the proof. The proof of part~b) is similar.
\end{proof}

Let $a$ and $h$ be paths as above. For $1\leq i\leq s$ denote
$v_i=a_i''$. We assume that the path $a$ is closed and primitive, $s\geq2$, $b'\neq b''$,
and $b',b''\in\{v_2,v_k\}$ for some $k\in\{1,3,4,\ldots,s\}$. Schematically this is
depicted as
$$
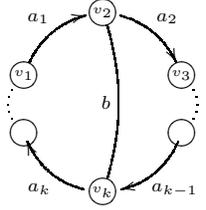
\vcenter{
\xymatrix@=.5cm{ %
&\vtx{v_2}  \ar@/_/@{<-}[ld]_{a_1}  \ar@/^/@{->}[rd]^{a_2}  \ar@/^/@{-}[ddd]_{b}& \\
\vtx{v_1}  \ar@/_/@{..}[d]  &&\vtx{v_3} \ar@/^/@{..}[d]\\
\vtx{\quad\,}  \ar@/_/@{<-}[rd]_{a_{k}}  &&\vtx{\quad\,} \ar@/^/@{->}[ld]^{a_{k-1}}\\
&\vtx{v_k}    & \\
}}.
$$

\begin{lemma}\label{lemma_L4}
Using the preceding notation we have $h\equiv a_1a_2 f_1\, a_1a_2 f_2$ for some
$f_1,f_2\in\path(\Q)$.
\end{lemma}
\begin{proof}
If $s=2$, then see part~b) of Lemma~\ref{lemma_Q_with 2ver}.

Suppose $b'=v_k$, $b''=v_2$, and $s\geq3$.  By Lemma~\ref{lemma_aaa},
\begin{eq}\label{eq_lemma_L4}
h\equiv ba_k a_{k+1}\cdots a_s f
\end{eq}
for some path $f$. We denote by $\QG$ the $h$-restriction of $\Q$ to the vertices
$v_1,v_2,v_3$ and consider $h$ as a path in $\QG$. Part~b) of Lemma~\ref{lemma_Q_with
3ver} together with~\Ref{eq_lemma_L4} concludes the proof. The case of $b'=v_2$ and
$b''=v_k$ is similar.
\end{proof}

\begin{lemma}\label{lemma_pA291}
Let $h$ be a closed path in a quiver $\Q$ and $h\not\equiv0$.
Then there exist pairwise different primitive closed paths $b_1,\ldots,b_r$,
$c_1,\ldots,c_t$ in $\Q$, where $r,t\geq0$, such that
$$\mdeg(h)=\sum_{i=1}^r \mdeg(b_i)+2\sum_{k=1}^t\mdeg(c_k);$$
and there are pairwise different arrows $x_1,\ldots,x_r$, $y_1,\ldots,y_t$,
$z_1,\ldots,z_t$ in $\Q$ satisfying
\begin{eq}\label{eq_pA291_1}
y_j,z_j\in\Arr{c_j}\text{ and }\deg_{y_j}(h)=\deg_{z_j}(h)=2,
\end{eq}
\vspace{-0.5cm}
\begin{eq}\label{eq_pA291_new0}
x_i\in\Arr{b_i}\text{ and }\deg_{x_i}(h)-2\sum_{k=1}^t \deg_{x_i}(c_k)=1 
\end{eq}
for any $1\leq i\leq r$, $1\leq j\leq t$.
\end{lemma}
\begin{proof}We assume that $\un{\de}=\mdeg(h)$.

If there is a $\un{\de}$-double path $a$ in $\Q$, then we define $c_1=a$. Let~\Ref{eq_pA291_1} be not valid for any $y_1,z_1\in\Arr{a}$ with $y_1\neq z_1$, i.e., there exists a $y\in\Arr{a}$ such that $\deg_{z}(h)\geq3$ for all $z\in\Arr{a}$ with $z\neq y$. Without loss of generality we can assume that $a=a_1\cdots a_s$, where $a_1,\ldots,a_s\in\Arr{\Q}$ are pairwise different and $a_1=y$. Then $h\equiv a_1a_2 f_{11} a_1a_2 f_{12}$ for some $f_{11},f_{12}\in\path(\Q)$ (see Lemma~\ref{lemma_L4}). Considering $a_1a_2$ as a new arrow, we can apply Lemma~\ref{lemma_L4} once again and obtain $h\equiv a_1a_2a_3 f_{21} a_1a_2a_3 f_{22}$ for some $f_{21},f_{22}\in\path(\Q)$. Repeating this procedure we can see that $h\equiv a f_{s-1,1} a f_{s-1,2}$ for some $f_{s-1,1},f_{s-1,2}\in\path(\Q)$. Since $a$ is a closed path, $h\equiv a^2 f_{s-1,1} f_{s-1,2}\equiv0$; a contradiction. Thus there are arrows $y_1,z_1$ satisfying the required conditions.  

We diminish $\un{\de}$ by $2\mdeg(a)$ and repeat the reasoning to obtain $c_2$ with $y_2,z_2$ and so on. Finally, we obtain $c_j,y_j,z_j$ for all $1\leq j\leq t$ $(t\geq0)$ such that the required conditions are valid and there is no $\un{\de}$-double path in $\Q$ for $\un{\de}=\mdeg(h)-2\sum_{k=1}^t\mdeg(c_k)$.

Assume $\un{\de}\neq0$. Since
\begin{eq}\label{eq_pA291_new1}
\sum_{a'=v}\de_a = \sum_{a''=v}\de_a,
\end{eq}%
there is a $\un{\de}$-single path $b$ in $\Q$ with $x\in\Arr{b}$ satisfying $\de_{x}=1$. We set $b_1=b$, $x_1=x$ and diminish $\un{\de}$ by $\mdeg(b)$. Repeating this procedure we obtain the required $b_1,\ldots,b_r$ together with $x_1,\ldots,x_r$ $(r\geq0)$.
\end{proof}

\begin{cor}\label{cor_new}
Suppose $\Q$ is a quiver with $d$ arrows and $m(\Q)=m$. Let $h$ be a closed path in $\Q$ and $h\not\equiv0$. Then $\deg(h)\leq md$.
\end{cor}
\begin{proof} We use notations from the formulation of Lemma~\ref{lemma_pA291}. Since $\deg(b_i)\leq m$ for $1\leq i \leq r$ and $\deg(c_k)\leq m$ for $1\leq k \leq t$,  we have 
$\deg(h)\leq m(r+2t)$. Moreover, there are $r+2t$ pairwise different arrows in $\Q$. Thus $r+2t\leq d$. 
\end{proof}

\section{Examples}\label{section_example}

Suppose $\Q$ is a strongly connected quiver. The {\it support} of a non-zero vector $\un{\de}\in\NN^{\#\Arr{\Q}}$ with respect to $\Q$ is the subquiver
$\Q_{\un{\de}}$ of $\Q$ such that $\Arr{\Q_{\un{\de}}}=\{a\in\Arr{\Q}\,|\,\de_a\geq1\}$ and $\Ver{\Q_{\un{\de}}}=\{a',a''\,|\,a\in \Arr{\Q_{\un{\de}}}\}$. We will apply the following remark together with Lemma~\ref{lemma_mdeg} to construct indecomposable invariants.

\begin{lemma}\label{lemma_inclusions_new} 
Let $\Char(K)=2$ and $h$ be a closed path in $\Q$. If for any $\mdeg(h)$-double path $a$ we have that the support of $\mdeg(h)-2\mdeg(a)$ is not strongly connected (and is not empty), then $h\not\equiv0$. 
\end{lemma}
\begin{proof} 
If $h$ satisfies the condition of the lemma and $h\equiv0$, then 
$h\equiv a^2f$ for some closed paths $a,f$. Thus the support of $\mdeg(h)-2\mdeg(a)=\mdeg(f)$ is strongly connected; a contradiction.
\end{proof}

%
%

\begin{lemma}\label{lemma_example}
Suppose $n\geq m\geq2$. Then for $d$ sufficiently large there is a quiver $\Q\in\QuiverNull(n,d,m)$ and a closed path $h$ in $\Q$ such that $h\not\equiv0$ and
\begin{enumerate}
\item[1)] $\deg(h) = md - (2nm - m^2 - m)$, if $\Char(K)=2$;

\item[2)] $\deg(h) =  3n$, if $\Char(K)\neq2$.
\end{enumerate}
\end{lemma}
\begin{proof} 
\textbf{1)} Suppose $\Char(K)=2$. For $d\geq 2n-m$ we consider the following strongly connected quiver $\Q\in \QuiverNull(n,d,m)$:
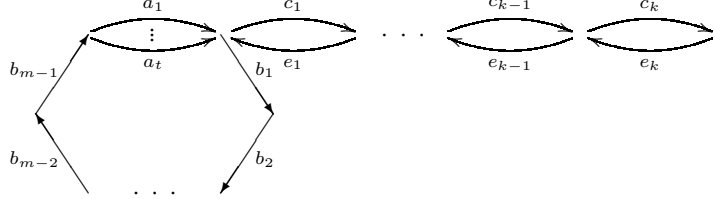
$$
\begin{picture}(0,80)
\put(-145,35){%
\put(0,0){\vector(2,3){20}}%
\put(20,-30){\vector(-2,3){20}}%
\put(70,30){\vector(2,-3){20}}%
\put(90,0){\vector(-2,-3){20}}%
\put(45,-30){%
\put(0,0){\circle*{1}}%
\put(7,0){\circle*{1}}%
\put(-7,0){\circle*{1}}}%
\put(15,30){\xymatrix@C=1.65cm@R=1cm{ %
\ar@/^/@{->}[r]^{a_1} \ar@/_/@{->}[r]_{a_t}&
\ar@/^/@{->}[r]^{c_1} \ar@/_/@{<-}[r]_{e_1}&\\
}}%
\put(150,30){\xymatrix@C=1.65cm@R=1cm{ %
\ar@/^/@{->}[r]^{c_{k-1}} \ar@/_/@{<-}[r]_{e_{k-1}}&
\ar@/^/@{->}[r]^{c_k} \ar@/_/@{<-}[r]_{e_k}&\\
}}%
\put(139,30){
\put(0,0){\circle*{1}}%
\put(7,0){\circle*{1}}%
\put(-7,0){\circle*{1}}%
}%
\put(44,30){\put(0,0){\circle*{1}}\put(0,2){\circle*{1}}\put(0,-2){\circle*{1}}}%
\put(300,-1){,}%
\put(-10,15){$\scriptstyle b_{m-1}$}%
\put(-10,-19){$\scriptstyle b_{m-2}$}%
\put(83,15){$\scriptstyle b_1$}%
\put(83,-19){$\scriptstyle b_2$}%
}%
\end{picture}$$%
where $t=d+m+1-2n\geq1$ and $k=n-m\geq0$. Here we assume that if $k=0$, then there are not  arrows $c_1,e_1,\ldots,c_k,e_k$. We take $h= b a_1 \cdots b\, a_t$, where $b=b_1\cdots b_{m-1}$. By Lemma~\ref{lemma_inclusions_new}, $h\not\equiv0$. Obviously,  $\deg(h)$ satisfies the required equality.  

\smallskip
\textbf{2)} Suppose $\Char(K)\neq 2$. For $d\geq 3n$ we consider the quiver from part 1). Then we remove arrows $a_2,\ldots,a_t$ from it and add one loop to each of the vertices $a_1',c_1',\ldots,c_{k-1}'$; we also add two loops to each of the rest of vertices. The resulting quiver is denoted by $\QG$. Then we add $d-3n$ arbitrary arrows to construct the required quiver $\Q\in\Q(n,d,m)$. By Lemma~\ref{lemma_mdeg}, there is a closed path $h$ in $\Q$ of degree one in each of the arrows of $\QG$. Thus, $\deg(h)=3n$. Since $\deg_a(h)\leq1$ and $\deg_v(h)\leq3$ for all $a\in\Arr{\Q}$ and $v\in\Ver{\Q}$, it is not difficult to see that the definition of the equivalence $\equiv$ implies that $h\not\equiv0$. 
\end{proof}

\begin{example}\label{ex_1}
We assume that $\Char(K)=2$ and consider the quiver prom part~1) of the proof of Lemma~\ref{lemma_example}.  Denote $b=b_1\cdots b_{m-1}$, $c=c_1\cdots c_{k-1}$, and $e=e_{k-1}\cdots e_1$. We take $h= b a_1 \cdots b\, a_t c\, e\, c\, c_k e_k e$. By Lemma~\ref{lemma_inclusions_new}, $h\not\equiv0$.

We set $B=X_{b_{m-1}}\cdots X_{b_1}$, $C=X_{c_{k-1}}\cdots X_{c_1}$, and $E=X_{e_{1}}\cdots X_{e_{k-1}}$. Then the invariants 
$$\tr(X_{a_1} B \cdots X_{a_t}B)\text{ and }\tr(X_{a_1} B \cdots X_{a_t} B\, E\, C\, E\, X_{e_k} X_{c_k} C)$$ 
are indecomposable by Lemma~\ref{lemma_reduction}. 
\end{example}


\bigskip
\noindent{\bf Acknowledgements.} This paper was written during author's visit to
Bielefeld University, sponsored by DAAD grant for young scientists (Forschungsstipendien
${\rm f\ddot{u}r}$ Doktoranden und Nachwuchswissenschaftler). The author is grateful for
this support. The author would like to thank Claus Michael Ringel for hospitality and   the anonymous referee for useful suggestions and comments. The paper has also been partially supported by RFFI 08-01-00067.



\begin{thebibliography}{99}
\bibitem{Amitsur_1980} Amitsur, S.A. (1980). On the characteristic polynomial of a sum of
matrices. {\it Linear and Multilinear Algebra} 8:177--182.

\bibitem{ADS_2006}Aslaksen, H., Drensky, V., Sadikova, L. (2006). Defining relations of invariants of two $3\times 3$ matrices. {\it J. Algebra} 298:41--57.

\bibitem{Drensky_survey_2007}Drensky, V. (2007). Computing with matrix invariants. {\it Math. Balkanica (N.S.)} 21(1-2):141--172. 

\bibitem{Drensky_Sadikova_4x4}Drensky, V., Sadikova, L. (2006). Generators of invariants of two $4\times 4$ matrices. {\it C. R. Acad. Bulgare Sci.} 59(5):477--484.

\bibitem{Domokos_1998}Domokos, M. (1998). Invariants of quivers and wreath products. {\it 
Comm. Algebra} 26:2807--2819.

\bibitem{Domokos_gen_2002}Domokos, M. (2002). Finite generating system of matrix invariants. {\it Math. Pannon.} 13(2):175--181.

\bibitem{DKZ_2002}Domokos, M., Kuzmin, S.G., Zubkov, A.N. (2002). Rings of
matrix invariants in positive characteristic. {\it J. Pure
Appl. Algebra} 176:61--80.  %

\bibitem{Donkin_1985}Donkin, S. (1985). {\it Rational representations of algebraic groups:
tensor products and filtrations}. Lecture Notes in Math., 1140. Berlin, Heidelberg, New York: Springer.

\bibitem{Donkin_1992a}Donkin, S. (1992). Invariants of several matrices. {\it Invent.
Math.} 110:389--401.

\bibitem{Donkin_1993}Donkin, S. (1993). On tilting modules for algebraic groups, {\it 
Math. Z.} 212:39--60.

\bibitem{Donkin_1994}Donkin, S. (1994). Polynomial invariants of representations of
quivers. {\it Comment. Math. Helvetici} 69:137--141.

\bibitem{Formanek_1987}Formanek, E. (1987). The invariants of $n\times n$ matrices. {\it 
Lecture Notes in Math.} 1278:18--43.

\bibitem{Formanek_1991}Formanek, E. (1991). {\it The polynomial identities and
invariants of $n\times n$ matrices}. Regional Conference series in Mathematics, 78. 
Providence, RI: American Math. Soc.

\bibitem{Gabriel_1972}Gabriel, P. (1972). Unzerlegbare Darstellungen I. {\it Manuscr.
Math.} 6:71--103.

\bibitem{Higman_1956}Higman, G. (1956). On a conjecture of Nagata. {\it 
Proc. Cambridge Philos. Soc.} 52:1--4.

\bibitem{Klein_2000}Klein, A.A. (2000). Bounds for indices of nilpotency and
nility. {\it Arch. Math. (Basel)} 76:6--10.

\bibitem{Le_Bruyn_Procesi_1990}Le Bruyn, L., Procesi, C. (1990). Semi-simple
representations of quivers. {\it Trans. Amer. Math. Soc.} 317:585--598.

\bibitem{Lopatin_Comm1}Lopatin, A.A. (2004). The algebra of invariants of $3\times 3$ matrices over a field of arbitrary characteristic. {\it Comm. Algebra} 32(7):2863--2883.

\bibitem{Lopatin_Comm2}Lopatin, A.A. (2005). Relatively free algebras with the identity $x^3=0$. {\it Comm. Algebra} 33(10):3583--3605.

\bibitem{Nakamoto_2002}Nakamoto, K. (2002). The structure of the invariant ring of two matrices of degree 3. {\it J. Pure Appl. Algebra} 166:125--148.

\bibitem{Procesi_1976}Procesi, C. (1976). The invariant theory of $n\times n$
matrices. {\it Adv. Math.} 19:306--381.

\bibitem{Procesi_1984}Procesi, C. (1984). Computing with $2\times 2$ matrices. {\it 
J. Algebra} 87:342--359.

\bibitem{Razmyslov_1974}Razmyslov, Yu.P. (1974). Trace identities of full matrix
algebras over a field of characteristic $0$. {\it Izv. Akad. Nauk SSSR Ser. Mat.} 38(4):723--756 (Russian).

\bibitem{Sibirskii_1968}Sibirskii, K.S. (1968). Algebraic invariants of a system of
matrices. {\it Sibirsk. Mat. Zh.} 9(1):152--164 (Russian).

\bibitem{Teranishi_1986}Teranishi, Y. (1986). The ring of invariants of matrices. {\it 
Nagoya Math. J.} 104:149--161.

\bibitem{Zubkov_Fund_Math_2001}Zubkov, A.N. (2001). The Razmyslov-Procesi theorem for
quiver representations, {\it Fundam. Prikl. Mat.} 7(2):387--421 (Russian).


\end{thebibliography}
\end{document}